\documentclass[12pt,oneside, reqno]{amsart}
\usepackage{amsmath}

\usepackage{tikz-cd}
\usepackage{graphicx}
\usepackage{hyperref}
\usepackage{amsfonts}
\usepackage{amssymb}

 \newtheorem{thm}{Theorem}[subsection]

 \newtheorem{lem}[thm]{Lemma}
 \newtheorem{prop}[thm]{Proposition}

\theoremstyle{definition}
 \newtheorem{defn}[thm]{Definition}
  \newtheorem{exm}[thm]{Example}
  
 \theoremstyle{remark}
 \newtheorem{rem}[thm]{Remark}
 
\numberwithin{equation}{section}

\newcommand{\R}{\mathbb{R}}
\newcommand{\D}{\mathcal{D}}
\newcommand{\F}{\mathcal{F}}

\newcommand{\dd}{\mathrm{d}}

\newcommand{\im}{\mathop{\mathrm{im}}}

\newcommand{\Ho}{\mathrm{H}}
\newcommand{\id}{\mathrm{id}}
\newcommand{\C}{\mathcal{C}}

\begin{document}
\title []{Mayer-Vietoris sequence for generating families in diffeological spaces}
\thanks{The authors were partially supported by  Xunta de Galicia ED431C 2019/10 with FEDER funds. The first author was partially supported by MINECO-Spain research project PID2020-MTM-114474GB-I00}
\author{Enrique Mac\'ias-Virg\'os, Reihaneh Mehrabi}
\address{Department of Mathematics, University of Santiago de Compostela, 15782-Spain}
\email{quique.macias@usc.es, reihaneh.mehrabi@rai.usc.es}


\maketitle 

\begin{abstract}We prove a version of the Mayer-Vietoris sequence for De Rham differential forms in diffeological spaces. It is based on the notion of a generating family instead of that of a covering by open subsets.
\end{abstract}

{MSC 2020: 14F40
58A40} 

\section*{Introduction}
Diffeological spaces are a generalization of differentiable manifolds which provides a  unified framework for non-classical objects in differential geometry  like quotients of manifolds, spaces of leaves of foliations, spaces of smooth functions or groups of diffeomorhisms. Roughly speaking, a diffeology on an  set $X$ specifies
which of the maps from a domain in $\R^m$, $m\geq 0$,  into $X$ (called {\em plots})  are smooth. In this setting one can define many constructions and invariants analogous to the classical ones, but with a much larger scope.
In particular, the differential graded algebra $\Omega(X)$ of  differential forms and the De Rham
cohomology can be generalized to this context.
The standard reference for the subject of diffeological spaces is Iglesias-Zemmour's book \cite{PATRICK}.

In the classical De Rham theory, the Mayer-Vietoris sequence is considered the main basic tool for computing the cohomology groups of a manifold that can be written as the union of two open subspaces.
In diffeological spaces, versions of this Mayer-Vietoris result have been given. For instance, Iwase and Izumida \cite{IWASE-IZUMIDA} introduced a new version of differential forms (cubical differential forms) in order to obtain partitions of unity and an exact Mayer-Vietoris sequence. Also, Haraguchi \cite{HARAGUCHI} considered the Mayer-Vietoris sequence for general differential forms and gave sufficient conditions for the existence of partitions of unity. In \cite{KURIBAYASHI},  Kuribayashi developed cohomological methods in diffeological spaces and used a Mayer-Vietoris argument for different cochain complexes.

In all these cases, the classical theory is imitated by considering two open sets of the so-called $\D$-topology. Those are the sets $W\subseteq X$ such that $\alpha^{-1}(W)$ is open in $U$ for any plot $\alpha\colon U \to X$ of the diffeology. However, it is our opinion that introducing the $\D$-topology fits poorly with the general philosophy of diffeological spaces, where topology is often irrelevant and differential objects are defined intrinsically: we can give a space a differentiable structure without first giving it a topology. See \cite{STACEY} for a discussion on this issue. In other words, smoothness is a property that is not based on continuity: for instance, the $\D$-topology is determined by the smooth curves, while diffeologies are not  \cite[Theorem 3.7]{C-S-W}.

There is another notion  more intimately linked to the foundations of the theory, namely that of a {\em generating family} (Definition \ref{GENER}). Many diffeologies are constructed by only giving a generating family of plots, often much smaller than the whole diffeology.  This mode of construction of diffeologies is very
useful because it  reduces the study of a diffeological space to
a subset of its plots.

Then, instead
of a covering $\{U,V\}$ of $X$ by two open sets, which would give the classical sequence
$$0\to \Omega(X) \to \Omega(U)\oplus \Omega(V) \to \Omega(U\cap V) \to 0,$$ our version of the Mayer-Vietoris result is based on a generating family $\{\alpha,\beta\}$ formed by two plots $\alpha,\beta$. 
We get an exact sequence
\begin{equation*}
0\to \Omega(X) \to \Omega(\alpha)\oplus \Omega(\beta) \to \Omega(P),
\end{equation*}
where $\Omega(\alpha)$  (resp. $\Omega(\beta)$) is  the subcomplex of $\Omega(U)$  of $\alpha$-horizontal (resp. $\beta$-horizontal) forms  (see \ref{HORIZ}) and $P$ is the pullback of two plots (see \ref{PULLBACKS}).  This idea is natural, because $U\cap V$ can be viewed as the pullback of the inclusions $U,V \subseteq X$.

 Our objective is to understand how the De Rham complex of $X$ is determined by $\alpha$ and $\beta$. This is achieved in our main Theorem \ref{MAIN}. Our procedure  is close to the diffeological gluing procedure in \cite{PERVOVA}, which can be considered as a form of Mayer-Vietoris construction.

Note that we do not obtain a long exact sequence in cohomology, which would need  the existence of partitions of unity, a practically impossible objective in such a general context. This is coherent with the fact that diffeological spaces is a wide category that allows abstract constructions but where the results that depend on the local structure  need additional hypothesis \cite{B-H}.

The contents of the paper are as follows. In Section \ref{DIFF-SPACE} we present the basic definitions that we need along the paper, including that of a generating family and the pullback of two plots. In Section \ref{DE-RHAM} we explain Cartan calculus and De Rham cohomology in diffeological spaces, and we introduce the notion of horizontal forms. In Section \ref{MAYER-VIET} we state our main Theorem \ref{MAIN}, which is the Mayer-Vietoris sequence associated to a generating family by two plots. Finally, in Section \ref{AN-EXAMPLE}, we compute the cohomology of an illustrative example.

\section{Diffeological spaces}\label{DIFF-SPACE}
An open subset $U\subseteq \R^m$ of some Euclidean space, $m\geq 0$, is called an {\em $m$-domain}. A differentiable $\C^\infty$ map $h\colon V\subseteq \R^n \to U\subseteq \R^m$ between domains will be called a {\em change of coordinates}. Let $X$ be a set. Any set map $\alpha\colon U\subseteq\R^m \to X$ defined on an $m$-domain is called a {\em parametrization} on $X$. 

\subsection{Diffeology}
A {\em diffeology} (of class $\C^\infty$) on the set $X$ is a family of parametrizations satisfying the following axioms. 
\begin{defn}
Let $X$ be a set. A diffeology on X is any family $\D$ of parametrizations on $X$ such that 
\begin{enumerate}
\item
any constant parametrization on $X$ belongs to $\D$; 
\item
 if $\alpha\colon U\subseteq\R^{n}\rightarrow X$ is a parametrization that locally belongs   to $\D$ (i.e. for every $x\in U$ there exists a neighbourhood $U_x$ such that $\alpha_{\vert U_x}\in \D$), then $\alpha\in D$;
\item
 if $\alpha \in \D$ and $h\colon V\to U$ is a change of coordinates, then $\alpha \circ h\in \D$.
\end{enumerate}
\end{defn}

A set endowed with a diffeology is called a {\em diffeological space}.
 The parametrizations of the diffeology $\D$ of $X$ are called {\em plots}. For the sake of clarity, sometimes we shall make explicit the domain of a plot $\alpha\colon U \to X$ by denoting it as $(U,\alpha)$.

\begin{exm}If $U\subseteq \R^m$ is an $m$-domain, the $\C^\infty$ changes of coordinates $V \to U$ with codomain $U$ form its {\em usual} diffeology. More generally, if $M$ is a finite dimensional differentiable manifold, the $\C^\infty$ maps $V\subseteq \R^m \to M$, with domain any $m$-domain, form the {\em usual} or {\em standard} diffeology on $M$. In this case, all the diffeological objects that we shall define correspond to the usual ones.
\end{exm}

\subsection{Differentiable maps}

\begin{defn}
A map $(X,\D)\overset{f}\rightarrow (X,\D')$ between diffeological spaces is {\em differentiable} if $f\circ\alpha\in \D'$ for all $\alpha \in \D$.
\end{defn}

The composition of differentiable maps is differentiable. The plots $U\to X$ of $\D$ are differentiable for the usual diffeology on $U$.

\subsection{Subspaces} A nice characteristic of diffeological spaces is that any subspace inherits a diffeology, making the theory much more clean than the usual theory of differentiable manifolds and submanifolds.

\begin{defn}If $(X,\D_X)$ is a diffeological space and $Y\subseteq X$ is any subset, we define the {\em induced diffeology} or {\em subspace diffeology} $\D_Y$ on $Y$ as the family of plots $(U,\alpha)$ in $\D_X$ whose image is contained in $Y$.
\end{defn}

With this diffeology, not only the inclusion $\iota_Y\colon Y \hookrightarrow X$ is differentiable, but a map $F\colon Z \to Y$ is differentiable if and only if the composition $\iota_Y\circ F \colon Z \to X$ is differentiable.

\subsection{Comparing diffeologies}
We shall avoid the words {\em finer} and {\em coarser} when comparing diffeologies on the same set $X$. If $\D \subseteq \D'$ we say that $\D'$ is {\em larger} than $\D$, equivalently $\D$ is {\em smaller} than $\D'$.

The smallest diffeology on $X$ is the {\em discrete} diffeology, formed by all the {locally constant} parametrizations on $X$.

The largest diffeology on $X$ is formed by all possible parametrizations. It is called the {\em indiscrete} or {\em coarse} diffeology.

\begin{prop}
The intersection of an arbitrary collection of diffeologies on $X$ is a diffeology on $X$.
\end{prop}

\subsection{Generating families}
The diffeologies can be built starting with a generating family. We refer to \cite[Art. 1.66 and fol.]{PATRICK} for the proofs of the basic results that we will need here.

\begin{defn}\label{GENER}
Let $\F=\{\alpha_j\}$ be an arbitrary family of parametrizations $\alpha_j\colon U_j\to X$ and let $\D=\langle\mathcal{F}\rangle$ be the intersection of all the diffeologies on $X$ containing $\F$. We say that $\mathcal{F}$ is a generating family for $\D$.
\end{defn}

The diffeology $\D=\langle F \rangle$ generated by the family $\F$ is then the smallest diffeology containing it.  The intersection is not void because $\F$ is always contained in the indiscrete diffeology. It is always possible to add to a generating family any family of parameterizations that generate the discrete diffeology, without altering the generated diffeology; consequently, we shall always assume that the family $\F$ contains all the locally constant parametrizations. Then $\langle \F\rangle$ can be characterized by the following property, that we will use heavily.

\begin{thm}\label{LOCALLY}
The parametrization $(U,\alpha)$ is a plot in $\langle\F\rangle$ if and only if for all $x\in U$ there exists a neighborhood $V=U_x$ of $x$ in $U$ such that  there exist a plot $f\colon W\to X$ in the family $\F$ and a change of coordinates $h\colon V \to W$  such that  $\alpha_{\vert V}=f\circ h$.
\end{thm}

\begin{exm}The usual diffeology on $\R^m$ is generated by the plot $\id\colon \R^m \to \R^m$.
\end{exm}

\subsection{Pullbacks}\label{PULLBACKS}
The category of diffeological spaces and differentiable maps allows any kind of categorical limits \cite{B-H}.  In particular, we shall need the following construction of the pullback (or {\em fibered product}) of two plots.

Let $(U,\alpha)$ and $(V,\beta)$ be two plots on the diffeological space $(X,\D)$.
We consider the subspace $P$ (endowed with the induced diffeology) of the product
$U\times V$ (endowed with the usual diffeology) given by
$$P=\{(u,v)\in U \times V \colon \alpha(u)=\beta(v)\};$$
also we consider the differentiable maps
$$p_U\colon P \to U, \quad p_U(u,v)=u,$$
and
$$p_V\colon P \to V, \quad p_V(u,v)=v,$$
induced by the projections. They verify
$$\alpha\circ p_U=\beta\circ p_V.$$

\begin{prop}\label{PULLBACKP} $P$ is a pullback in the category of diffeological spaces. That means that it verifies the following universal property: given any diffeological space $Y$ and two differentiable maps $a\colon Y \to U$ and $b\colon Y \to V$ such that $\alpha\circ a = \beta \circ b$, then there exists a unique differentiable map $F\colon Y \to P$ such that $p_U\circ F=a$ and $p_V\circ F = b$,

\begin{center}
\begin{tikzcd}
Y\arrow[ddr,bend right,"a"']\ar[drr,bend left,"b"]\ar[dr,dotted,"F"]&&\\
&P\ar{d}{p_U}\ar{r}{p_V}&V\subseteq\mathbb{R}^n \arrow[d,"\beta"]  \\
&U\subseteq\mathbb{R}^m  \ar{r}{\alpha}&X 
\end{tikzcd}
\end{center} \end{prop}

\section{De Rham
cohomology}\label{DE-RHAM}
Exterior differential forms and De Rham
cohomology can be generalized to the  context of diffeological spaces.   We refer to \cite{BOTT-TU,MADSEN} and \cite[Art. 6.28 and fol.]{PATRICK} for the basics on Cartan-De Rham Calculus.

\subsection{Differential forms}
A differential form on a diffeological space is defined as
the family of its pullbacks by the plots of the diffeology, which are usual differential forms in $m$-domains.

\begin{defn}
Let $(X,\D)$ be a diffeological space. A differential $k$-form on $X$, $k\geq 0$, is any collection $w=\{w_\alpha\}$, where for each plot $\alpha\colon U\to X$ in $\D$ we have a usual differential form $w_\alpha\in \Omega^k(U)$. Moreover, it must verify the following compatibility condition:   for any $\C^\infty$ change of coordinates $h\colon V\to U$ we must have 
$$w_{\alpha \circ h}=h^*w_\alpha,$$
where $h^*\omega_\alpha\in \Omega^k(V)$ is the usual pullback of a differential form.
\end{defn}

We shall denote by $\Omega^k(X,\D)$ the vector space of $k$-forms on $(X,\D)$. It is another nice feature of diffeological theory that this set of forms can be endowed with a diffeology in a natural way (\cite[Art. 6.29]{PATRICK}.

\begin{rem}\label{VERTICAL}The form $\omega_\alpha$ is not an arbitrary form on $U$. Notice that if $h,h'\colon V \to U$ are two changes of coordinates such that $\alpha\circ h = \alpha\circ h'$, then
$$h^*\omega_\alpha =\omega_{\alpha\circ h }=\omega_{\alpha\circ h'}=(h')^*\omega.$$
We shall need this property later (see Section \ref{HORIZ}).
\end{rem}

\begin{exm}
If $M$ is a differentiable manifold, the diffeological forms for the usual diffeology are the usual differential forms \cite{PATRICK,H-M-S}
\end{exm}

\begin{exm}\label{ONEPOINT}
Let $X=\{*\}$ be a one-point set (with the discrete diffeology). Then $\Omega^0(\{*\})=\R$ and $\Omega^k(\{*\})=0$ for $k\geq 1$.
\end{exm}

\begin{prop}[{\cite[Art. 6.31]{PATRICK}}]The zero forms on a diffeological space $(X,\D)$ are identified with the differentiable maps $X \to \R$, for the usual diffeology on $\R$.
\end{prop}
\begin{proof}
If $f\colon X \to \R$ is a differentiable map, we define the form $\omega(f)\in \Omega^0(X)$ as
$$\omega(f)_\alpha=f\circ \alpha\in \Omega^0(U)=\C^\infty(U,\R)$$
for any plot $(U,\alpha)$ on $X$.

Conversely, if $\omega\in \Omega^0(X)$, we define the function
$$f(\omega)\colon X \to \R, \quad f(\omega)(x)=\omega_{c_x},$$
where for each $x\in X$ we take the constant plot $c_x\colon \R^0=\{0\}\to X$ with $c_x(0)=x$. 

The function $f=f(\omega)$ is differentiable because for any plot $(U,\alpha)$ on $X$, we have
$$(f\circ \alpha)(u)=\omega_{c_{\alpha(u)}}=\omega_{\alpha\circ c_u}=c_u^*(\omega_\alpha)=\omega_\alpha(u),$$
where $c_u\colon \R^0\to U$ is the constant map $u$, and $\omega_\alpha\in\C^\infty(U,\R)$.\end{proof}

%


\begin{defn}
If  $F\colon (X,\D)\rightarrow (X',\D')$ is a differentiable map between diffeological spaces and $\omega\in \Omega^k(X',\D')$ is a $k$-form on $X'$  then we define the pull-back $F^*\omega\in \Omega^{k}(X,\D)$ as
 $$(F^*\omega)_{\alpha}:=\omega_{F\circ \alpha}$$
 for any plot $\alpha$ on $X$.
\end{defn}

Notice that if $\omega\in \Omega^k(X,\D)$, then $\omega_\alpha=\alpha^*\omega$ for any plot $\alpha\in\D$.


\begin{defn}The exterior derivative $\dd_X\colon \Omega^k(X,\D) \to \Omega^{k+1}(X,\D)$
is defined, for any plot $(U,\alpha)\in\D$,  by
$$(\dd_X\omega)_\alpha=\dd_U(\omega_\alpha),$$
where $\dd_U\colon \Omega^k(U) \to \Omega^{k+1}(U)$ is the usual exterior derivative.
\end{defn}

The most importan property of exterior differential is that $\dd_X\circ \dd_X=0$, as follows immediately from the  analogous property for $m$-domains.

The pullback behaves well with the usual constructions in Cartan calculus, among them the exterior differential: if $F\colon X \to X'$ is a differentiable map, then
$$\dd_{X}(F^*\omega)=F^*(\dd_{X'}\omega).$$

\subsection{De Rham cohomology}
We already have all the tools necessary to define a De Rham cohomology theory on diffeological spaces. 

A $k$ form $\omega\in \Omega^k(X)$ is {\em closed} if $\dd_X\omega=0$. It is {\em exact} if $\omega=\dd_X\mu$ for some $(k-1)$-form $\mu$.
Every exact form is closed. 
In other words, we have a {\em De Rham complex}
$$0\to \Omega^0(X,\D) \xrightarrow{d_0} \Omega^1(X,\D)\xrightarrow{d_1} \cdots \to \Omega^k(X,\D)\xrightarrow{d_k} \Omega^{k+1}(X,\D) \to \cdots$$
such that $\im \dd_{k-1}\subseteq \ker d_k$.
We define the De Rham's $k$-cohomology group of $(X,\D)$ as the quotient of vector spaces
$$\Ho^k(X,\D) = \ker d_k / \im d_{k-1}.$$
It measures the failure of closed forms to be exact.
\begin{exm}It was proved in \cite{H-M-S} that the cohomology of the leaf space $M/\F$ of a foliated manifold $(M,\F)$, endowed with the quotient diffeology, equals the so-called basic cohomology of the foliation.
\end{exm}
Each differentiable map $F\colon (X,\D) \to (X',\D')$ induces a morphism
$$F^*\colon \Ho^k(X',\D') \to \Ho(X,\D), \quad F^*([\omega])=[F^*\omega],$$
and the usual properties hold,
$$(G\circ F)^*=F^*\circ G^*, \quad \id_X^*=\id_{\Ho(X)}.$$

\subsection{Horizontal forms}\label{HORIZ}
\begin{defn}\label{VERTEQ}
In view of Remark \ref{VERTICAL}, if $(U,\alpha)$ is a plot on $X$ we  define the subspace $\Omega^k(\alpha)\subseteq \Omega^k(U)$ of {\em $\alpha$-horizontal  forms}  as
those forms $\mu\in \Omega^k(U)$ verifying: if $h,h'\colon V \to U$ are changes of coordinates such that $\alpha\circ h = \alpha\circ h'$, then
$ h^*\mu = (h')^*\mu$.
\end{defn}

They form  a subcomplex of $\Omega(U)$ because if $\omega\in \Omega^k(\alpha)$ and $\alpha\circ h=\alpha\circ h'$ then
$$h^*(\dd_U\omega)=\dd_Vh^*\omega=\dd_V(h')^*\omega=(h')^*\dd_U\omega,$$
hence $\dd_U\omega\in\Omega^{k+1}(\alpha)$.



\section{The Mayer-Vietoris sequence}\label{MAYER-VIET}
We shall consider a generating family of $(X,\D)$ formed by {\em two} plots $(U,\alpha)$ and $(V,\beta)$. We want to understand to what extent the cohomology $\Ho(X,\D)$ is
determined by these two plots. 

\subsection{Restriction and difference morphisms}
Let $(P,p_U,p_V)$ be the pullback defined in Section \ref{PULLBACKS} and Proposition \ref{PULLBACKP}:
\begin{equation}\label{DIAG}
\begin{tikzcd}
P\ar{d}{p_U}\ar{r}{p_V}&V\subseteq\mathbb{R}^n \arrow[d,"\beta"]  \\
U\subseteq\mathbb{R}^m  \ar{r}{\alpha}&X 
\end{tikzcd}
\end{equation}

We have the sequence of complexes
\begin{equation}\label{MY}
0\to  \Omega^k(X) \xrightarrow{r} \Omega^k(\alpha)\oplus\Omega^k(\beta) \xrightarrow{\delta}  \Omega(P)
\end{equation}
where we call
$r$ a {\em restriction morphism} and $\delta$ the {\em difference morphism}.  
They are defined  by
 $$r(\omega)=(\alpha^*\omega, \beta^*\omega)$$
 and
$$\delta(\mu,\nu)= p^*_U \mu-p^*_V\nu.$$
The complexes $\Omega(\alpha)$ and $\Omega(\beta)$ of horizontal forms were considered in Section \ref{HORIZ}.
\begin{lem}$r$ and $\delta$ are well defined morphisms of complexes.
\end{lem}
\begin{proof}
 If $\omega\in \Omega^k(X)$, and $h,h' \colon V \to U$ are two changes of coordinates such that $\alpha\circ h = \alpha\circ h'$, then
$$h^*(\alpha^*\omega)=(\alpha\circ h)^*\omega=(h')^*(\alpha^*\omega).$$
Then $\alpha^*\omega$ is a horizontal form.  Analogously for $\beta$. This proves that $r$ is well defined.

 The morphism $r$ commutes with the differentials, $(\dd_U\oplus \dd_V)\circ r=r\circ \dd_{X}$, because
\begin{align*}
(\dd_U\oplus \dd_V)r(\omega)=&(\dd_U\oplus \dd_V)(\alpha^*\omega,\beta^*\omega)
=
(\dd_U\alpha^*\omega,\dd_V\beta^*\omega)\\
=&(\alpha^*\dd_X\omega,\beta^*\dd_X\omega)=r(\dd_X(\omega).
\end{align*}
Finally,
\begin{align*}
(\dd_P\circ \delta)(\mu,\nu)=&\dd_P(p^*_U\mu-p^*_V\nu)= d_Pp^*_U\mu- d_Pp^*_V\nu
\\
=&p^*_U\dd_U\mu- p^*_Vd_V\nu=\delta   (d_U\mu, d_V\nu).
\end{align*}
This proves
$d_P\circ \delta=\delta \circ (d_U\oplus d_V)$.
\end{proof}

\begin{thm}\label{MAIN}The sequence \eqref{MY} is exact.
\end{thm}

\begin{proof}

 We shall make the proof in several steps.
 
Step 1. The morphism $r$ is injective.

Let $\omega\in \Omega^k(X)$ such that $r(\omega)=0$. We must show that $w_\gamma=0$ for any plot $\gamma \colon W \to X$. Since $\{\alpha, \beta\}$ is a generating family, for each  $p\in W$ there exists a neighborhood $W$ of $p$ such that either $\gamma_{\vert W}$ is constant, in which case  $\gamma$ factors through $c_0\colon W \to \R^0=\{*\}$ and we have $\omega_\gamma=c_0^*0=0$ on $W$, or there is a change of coordinates $h\colon W\to  U$ with $\gamma_{\vert W}=\alpha \circ h$ (analogously $\gamma_{\vert W}=\beta\circ h'$ for some $h'\colon W \to V$). Then, on $W$ we have
\begin{align*}
\omega_{\gamma}=\gamma^*\omega=(\alpha\circ h)^*\omega=h^*\alpha^*\omega=0.
\end{align*}
This proves that $\omega=0$.

Step 2. $ \im r\subseteq \ker \delta$.

 We show $\delta \circ r=0$:
\begin{align*}
p_U^*\alpha^*w-p_V^*\beta^*\omega=(\alpha \circ p_U)^*\omega-(\beta\circ p_V)^*\omega=0.
\end{align*}
because Diagram \ref{DIAG} is commutative.

Step 3.  $\ker \delta \subseteq \im r$.

Let $(\mu,\nu)\in \Omega^k(\alpha)\oplus\Omega^k(\beta)$ such that
$p_U^*\mu=p_V^*\nu$. We need to define a form $\omega\in \Omega^k(X)$ such that
$\alpha^*\omega= \mu$ and $\beta^*\omega=\nu$. For that, we need to define $\omega_\gamma\in \Omega^k(W)$ for any plot $(W,\gamma)$ on $X$.

Obviously, for $\gamma=\alpha$ we take $\omega_\alpha=\mu$. If $h\colon W\to U$ is a change of coordinates, we put $\omega_{\alpha\circ h}=h^*\mu$. This form  does not depend on $h$ due to the definition of horizontal form (Section \ref{HORIZ}), that is, if $\alpha\circ h=\alpha\circ h'$ then $h^*\mu=(h')^*\mu$. Analogously, for $\gamma=\beta$ we take $\omega_\beta=\nu$, and for a change of coordinates $h'\colon W \to V$ we put $\omega_{\beta\circ h'}=(h')^*\nu$. It may happen that $\alpha\circ h =\beta \circ h'$, for some changes of coordinates $h\colon W \to  U$ and $h'\colon W \to V$. In this case, by the pullback property of Proposition \ref{PULLBACKP}, there exists a differentiable map $F\colon W \to P$ such that
$p_U\circ F=h$ and $p_V\circ F=h'$. Hence
$$(h')^*\nu=(p_V\circ F)^*\nu=F^*p_V^*\nu=F^*p_U^*\mu=(p_U\circ F)^*\mu=h^*\mu.$$
This proves that the preceding definitions are consistent.

Finally, let $\gamma\colon W \to X$ be an arbitrary plot. For each point $p\in W$ there is a neighbourhood $W_p$ such that the restriction of $\gamma$ to $W_p$ either is constant, in which case we define $\omega_\gamma=0$ on the neighbourhood $W_p$, or $\gamma_{\vert W_p}=\alpha \circ h$ for some $h\colon W_p \to U$, in which case we define $\omega_\gamma=h^*\omega_\alpha$, or 
$\gamma_{\vert W_p}=\beta \circ h'$, in which case we state $\omega_\gamma$ to be $(h')^*\omega_\beta$ on $W_p$.

The compatibility condition can be checked easily. This proves that we have a well defined $k$ form $\omega$ on $X$ such that $r(\omega)=(\mu,\nu)$.\end{proof}

The latter Theorem shows that the complex $\Omega^k(X)$ is isomorphic to the subcomplex $\ker \delta$ of $\Omega(\alpha)\oplus\Omega(\beta)$.

Note that we do not prove the surjectiveness of $\delta$. This would require the existence of partitions of unity, a practically impossible objective in such a wide context.

\section{An Example}\label{AN-EXAMPLE}
We shall  compute the De Rham cohomology of the following Example as an application of our Theorem \ref{MAIN}. 

Let $X$ be the set obtained as the union of the two coordinate axis in $\R^2$, that is,
$$X=\{(x,y)\in \R^2\colon xy=0\}.$$
We take on $X$ the diffeology $\D$ generated by the two axis inclusions, that is, by the plots
$$\alpha\colon U=\R \to X, \quad \alpha(s)=(s,0),$$
and
$$\beta\colon V=\R \to X, \quad \beta(t)=(0,t).$$
The pullback of these two plots is $P=\{(0,0)\}$, because
$\alpha(s)=\beta(t)$ if and only if $s=t=0$. Hence Example \ref{ONEPOINT} applies.

\begin{rem}Notice that this diffeology $\D$ is strictly smaller than the subspace diffeology $\D'$ of $X$ as a subset of $\R^2$. In fact, the $1$-plot $\alpha\colon \R \to X\subseteq \R^2$ given by
the $\C^\infty$ map 
$$\alpha(t)=
\begin{cases}
(e^{1/t},0) \quad &\text{if\ } t<0,\\
(0,0) \quad &\text{if\ } t=0,\\
(0,e^{-1/t}) \quad &\text{if\ } t>0,
\end{cases}$$
does not belong to the diffeology, by Proposition \ref{LOCALLY}. Hence the computation of the $1$-differential forms for $(X,\D')$ (the so-called {\em cross})  in \cite{CROSS} does not apply here.
\end{rem}
We compute the horizontal forms for $\D$: since $p_1\alpha=\id_\R$, for the first projection $p_1\colon \R^2 \to \R$, the condition $\alpha\circ h = \alpha\circ h'$ for $h,h'\colon W \to U=\R$ means in fact $h=h'$. Hence $\Omega(\alpha)=\Omega(U)$. Analogously $\Omega(\beta)=\Omega(V)$.

Now, we know that $\Omega^k(X)\cong \im r=\ker\delta\subseteq \Omega^k(\alpha)\oplus\Omega^k(\beta)$. For $k\geq 2$ that means $\Omega^k(X)=0$. For $k=0$ we have that $(F,G)\in \Omega^0(\R)\oplus \Omega^0(\R)$ is a pair of functions $F,G\colon \R \to \R$, , and   being in the kernel of $\delta$ means $F(0)=G(0)$. For $k=1$, we have a pair of $1$-forms $(\mu,\nu)\in \Omega^1(\R)\oplus\Omega^1(\R)$, say
$\mu = f(t)\dd t$ and $\nu = g(t) \dd t$. The pair of functions $(f,g)$ is arbitrary because $p_U^*\mu=0=p_V^*\nu$ in $\Omega^1(P)=0$.

We consider the resulting De Rham complex of $X$, say,
$$0\to \Omega^0(X) \xrightarrow{\dd} \Omega^1(X) \to 0\to \cdots$$
which can be written as
$$0\to \C^\infty(\R)\oplus_0 \C^\infty(\R) \xrightarrow{\dd} \C^\infty(\R)\oplus \C^\infty(\R) \to 0 \to \cdots$$
where
$$\dd(F,G)=(F', G').$$
The kernel of $\dd$ is a pair of constant functions, which must be equal by the condition $F(0)=G(0)$, hence $$H^0(X,\D)=\R.$$
We compute the image of $\dd$. Since any $f\in \C^\infty(\R)$ is the derivative of some $F$, unique up to a constant, we can always write $f=\dd F$, $g= \dd G$ and we can assume that $F(0)=G(0)$, by taking $G-G(0)+F(0)$. Hence $\dd$ is surjective and
$$H^1(X,\D)=0.$$

\end{document}